\documentclass[10pt,twoside,reqno]{amsart}
\usepackage{amssymb}
\usepackage{multirow}
\usepackage{graphicx}
\usepackage[dvips]{epsfig}
\usepackage{latexsym}
\usepackage{amsmath}
\usepackage{amsthm}
\usepackage{amscd}
\usepackage{amssymb}
\input xy
\xyoption{all}

\calclayout

\numberwithin{equation}{section}
\makeatletter

\renewcommand{\@secnumfont}{\bfseries}

\renewcommand{\section}{\@startsection{section}{1}%
  {0mm}{.7\linespacing\@plus\linespacing}{.5\linespacing}
  {\normalfont\bfseries\centering}}

\newcommand{\bibsection}{\@startsection{section}{1}%
  {0mm}{.7\linespacing\@plus\linespacing}{.5\linespacing}
  {\normalfont\scshape\centering}}

\renewcommand{\@biblabel}[1]{#1.}

\newtheorem{thm}{\bf Theorem}[section]

\theoremstyle{theorem}

\theoremstyle{definition}

\numberwithin{equation}{section}

\begin{document}

\vspace{1.3cm}

\title {On type $2$ degenerate Bernoulli and Euler polynomials of complex variable}

\author{Taekyun Kim $^{1}$}
\address{$^{1}$ Department of Mathematics, Kwangwoon University, Seoul 139-701, Republic of Korea}
\email{tkkim@kw.ac.kr}

\author{Dae San Kim$^{2}$}
\address{$^{2}$ Department of Mathematics, Sogang University, Seoul 121-742, Republic of Korea}
\email{dskim@sogang.ac.kr}

\author{Lee-Chae Jang$^{3,*}$}
\address{$^{3}$ Graduate School of Education, Konkuk University, Seoul, 05029, Republic of Korea}
\email{Lcjang@konkuk.ac.kr}

\author{Han-Young Kim$^{4}$}
\address{$^{4}$ Department of Mathematics, Kwangwoon University, Seoul 139-701, Republic of Korea,}
\email{gksdud213@kw.ac.kr}

\thanks{$^{*}$ corresponding author}

\maketitle

\begin{abstract}

Recently, Masjed-Jamei-Beyki-Koepf studied the so called new type Euler polynomials without making use of Euler polynomials of complex variable. Here we study degenerate and type 2 versions of these
new type Euler polynomials, namely the type 2 degenerate cosine-Euler and type 2 degenerate sine-Euler
polynomials and also the corresponding ones for Bernoulli polynomials, namely the type 2 degenerate cosine-
Bernoulli and type 2 degenerate sine-Bernoulli polynomials by considering the degenerate Euler
and degenerate Bernoulli polynomials of complex variable and by treating the real and imaginary parts separately. We derive some explicit expressions for those new polynomials and some identities relating to them. Here we note that the idea of separating the real and imaginary parts separately gives an affirmative answer to the question asked by Hac\`ene Belbachir.
\end{abstract}

\pagestyle{myheadings}
\markboth{\centerline{\scriptsize T. Kim, D.S. Kim,L.C. Jang, H. Y. Kim }}
          {\centerline{\scriptsize On type $2$ degenerate Bernoulli and Euler polynomials of complex variable}}

\medskip
\section{\bf Introduction}
\medskip

As is known, the type $2$ Bernoulli polynomials $B_n(x)$, $(n\geq 0)$, and the type $2$ Euler
polynomials $E_n(x)$, $(n\geq 0)$, are respectively defined by 
\begin{equation}\begin{split}\label{eq01}
e^{xt} \frac{t}{2} csch \frac{t}{2} =\frac{t}{e^{\frac{t}{2}}-e^{-\frac{t}{2}}} e^{xt}
=\sum_{n=0}^\infty B_n(x) \frac{t^n}{n!},
\end{split}\end{equation}
and
\begin{equation}\begin{split}\label{eq02}
e^{xt} sech \frac{t}{2} =\frac{2}{e^{\frac{t}{2}}+e^{-\frac{t}{2}}}e^{xt}
=\sum_{n=0}^\infty E_n(x) \frac{t^n}{n!},  \quad \text{ (see \cite{ref04})}.
\end{split}\end{equation}
When $x=0$, $B_n=B_n(0)$ (or $E_n=E_n(0)$) are called the type $2$ Bernoulli (or type 2 Euler) numbers.

For $n\geq 0$, the central factorial numbers of the second kind are defined by the generating function to be
\begin{equation}\begin{split}\label{eq03}
\frac{1}{k!}\left(  e^{\frac{t}{2}} -e^{-\frac{t}{2}} \right)^k
=\sum_{n=k}^\infty T(n,k) \frac{t^n}{n!},  \quad \text{ (see \cite{ref02-1})}.
\end{split}\end{equation}
From \eqref{eq03}, we note that
\begin{equation}\begin{split}\label{eq04}
x^n = \sum_{k=0}^n T(n,k) x^{[k]} , \;(n\geq 0), \quad \text{ (see \cite{ref06})},
\end{split} \end{equation}
where $x^{[0]} =1$, $x^{[n]}=x\left(x+\frac{n}{2}-1\right)\left(x+\frac{n}{2}-2\right)
\cdots \left(x-\frac{n}{2}+1\right),\;(n\geq 1)$.
For $\lambda \in \mathbb{R}$, the degenerate exponential functions are defined as
\begin{equation}\begin{split}\label{eq05}
e^x_\lambda(t)=(1+\lambda t)^{\frac{x}{\lambda}} ,
e_\lambda (t)=e^1_\lambda (t)=(1+\lambda t)^{\frac{1}{\lambda}}.
\end{split}\end{equation}
By  \eqref{eq05}, we get
\begin{equation}\begin{split}\label{eq06}
e^x_\lambda(t)=\sum_{n=0}^\infty (x)_{n,\lambda} \frac{t^n}{n!} \quad \text{ (see \cite{ref07, ref08, ref09, ref10})},
\end{split}\end{equation}
where
\begin{equation}\begin{split}\label{eq07}
(x)_{0,\lambda}=1, \; (x)_{n,\lambda}=x(x-\lambda)\cdots(x-(n-1)\lambda), \;\;\;(n\geq 1).
\end{split}\end{equation}
In \cite{ref01, ref02}, Carlitz considered the degenerate Bernoulli polynomials given by
\begin{equation}\begin{split}\label{eq08}
\frac{t}{e_\lambda(t)-1} e_\lambda^x(t)
=\sum_{n=0}^\infty \beta_{n,\lambda}(x) \frac{t^n}{n!}.
\end{split}\end{equation}
When $x=0$, $\beta_{n,\lambda}=\beta_{n,\lambda}(0)$ are called the degenerate Bernoulli numbers.
In \cite{ref06}, Kim-Kim introduced the degenerate central factorial polynomials of the second kind which are given by
\begin{equation}\begin{split}\label{eq09}
\frac{1}{k!}\left(  e_\lambda^{\frac{1}{2}}(t)-  e_\lambda^{-\frac{1}{2}}(t)\right)^k
e_\lambda^x (t)= \sum_{n=k}^\infty T_\lambda(n,k |x) \frac{t^n}{n!},
\end{split}\end{equation}
where $k$ is a nonnegative integer.
When $x=0$, $T_\lambda(n,k)=T_\lambda(n,k|0)$ are called the degenerate
central factorial numbers of the second kind.

Recently, as a degenerate version of \eqref{eq01}, the type $2$ degenerate Bernoulli polynomials
are defined by 
\begin{equation}\begin{split}\label{eq10}
\frac{t}{e_\lambda^{\frac{1}{2}}(t)-  e_\lambda^{-\frac{1}{2}}(t)}
e_\lambda^x (t)= \sum_{n=0}^\infty B_{n,\lambda}(x)  \frac{t^n}{n!},\quad \text{ (see \cite{ref04})}.
\end{split}\end{equation}
When $x=0$, $B_{n,\lambda}=B_{n,\lambda}(0)$ are the type $2$ degenerate Bernoulli numbers.
By the same motivation as \eqref{eq10}, the type $2$ Euler polynomials are defined by
\begin{equation}\begin{split}\label{eq11}
\frac{2}{e_\lambda^{\frac{1}{2}}(t) + e_\lambda^{-\frac{1}{2}}(t)}
e_\lambda^x (t)= \sum_{n=0}^\infty E_{n,\lambda}(x)  \frac{t^n}{n!},\quad \text{ (see \cite{ref04})}.
\end{split}\end{equation}
When $x=0$, $E_{n,\lambda}=E_{n,\lambda}(0)$ are the type $2$ degenerate Euler numbers.\\

\indent Recently, several authors studied the degenerate Bernoulli and degenerate Euler numbers and polynomials (see \cite{ref01,ref02,ref03,ref04,ref05,ref06,ref07,ref08,ref09,ref10,ref11,ref12,ref13,ref14,ref15}). In addition,  Jeong-Kang-Rim introduced symmetry identities for Changhee polynomials of type two closely related to type 2 degenerate Euler polynomials (see \cite{ref04-1}), and Zhang and Lin obtained some interesting identities involving trigonometric functions and Bernoulli numbers (see \cite{ref15}).\\
\indent In \cite{ref05-1}, the authors considered the degenerate Bernoulli and degenerate Euler polynomials of complex variable. By treating the real and imaginary parts separately, they were able to introduce the degenerate cosine-Bernoulli polynomials, degenerate sine-Bernoulli polynomials, degenerate cosine-Euler polynomials and  degenerate sine-Euler polynomials, and derived some interesting results for them. \\
\indent In this paper, we study the type 2 degenerate Bernoulli and type 2 degenerate Euler polynomials of complex variable of which the latters are degenerate and type 2 versions of the new type Euler polynomials studied in \cite{ref12}. By treating the real and imaginary parts separately, the type 2 degenerate cosine-Bernoulli and type 2 degenerate sine-Bernoulli polynomials are introduced. We derive some explicit expressions for those polynomials and some identities relating to them. Moreover, the type 2 degenerate cosine-Euler and type 2 degenerate sine-Euler polynomials are investigated and analogous results to  the type 2 degenerate cosine-Bernoulli and type 2 degenerate sine-Bernoulli polynomials are obtained for them.

\medskip

\section{\bf  Type $2$ degenerate Bernoulli and Euler polynomials
of complex variable}
\medskip

From \eqref{eq10}, we define the type $2$ degenerate Bernoulli polynomials of complex
variable by 
\begin{equation}\begin{split}\label{eq12}
\frac{t}{e_\lambda^{\frac{1}{2}}(t) - e_\lambda^{-\frac{1}{2}}(t)}
e_\lambda^{x+iy} (t)= \sum_{n=0}^\infty B_{n,\lambda}(x+iy)  \frac{t^n}{n!},
\end{split}\end{equation}
and
\begin{equation}\begin{split}\label{eq13}
\frac{t}{e_\lambda^{\frac{1}{2}}(t) - e_\lambda^{-\frac{1}{2}}(t)}
e_\lambda^{x-iy} (t)= \sum_{n=0}^\infty B_{n,\lambda}(x-iy)  \frac{t^n}{n!},
\end{split}\end{equation}
where $i=\sqrt{-1}$. As is known, the degenerate cosine  and
sine functions are defined by
\begin{equation}\begin{split}\label{eq14}
cos_\lambda^{(y)} (t)=cos\left(  \frac{y}{\lambda} log(1+\lambda t) \right),
\end{split}\end{equation}
and
\begin{equation}\begin{split}\label{eq15}
sin_\lambda^{(y)} (t)=sin \left(  \frac{y}{\lambda} log(1+\lambda t) \right),\quad \text{ (see \cite{ref05-1})}.
\end{split}\end{equation}
Note that $\lim_{\lambda \rightarrow 0} cos_\lambda^{(y)}(t)=cos yt$,
$\lim_{\lambda \rightarrow 0} sin_\lambda^{(y)}(t)= sin yt$.
From \eqref{eq12} and \eqref{eq13}, we can derive the following equations.
\begin{equation}\begin{split}\label{eq16}
\sum_{n=0}^\infty \left(\frac{B_{n,\lambda}(x+iy)+B_{n,\lambda}(x-iy) }{2}  \right)\frac{t^n}{n!}
=\frac{t}{e_\lambda^{\frac{1}{2}}(t)-e_\lambda^{-\frac{1}{2}}(t)} e_\lambda^x(t) cos_\lambda^{(y)}(t),
\end{split}\end{equation}
and
\begin{equation}\begin{split}\label{eq17}
\sum_{n=0}^\infty \left(\frac{B_{n,\lambda}(x+iy)-B_{n,\lambda}(x-iy) }{2i}  \right)\frac{t^n}{n!}
=\frac{t}{e_\lambda^{\frac{1}{2}}(t)-e_\lambda^{-\frac{1}{2}}(t)} e_\lambda^x(t) sin_\lambda^{(y)}(t).
\end{split}\end{equation}
Now, we define the type $2$ degenerate cosine-Bernoulli and sine-Bernoulli polynomials 
by the generating functions as
\begin{equation}\begin{split}\label{eq18}
\frac{t}{e_\lambda^{\frac{1}{2}}(t)-e_\lambda^{-\frac{1}{2}}(t)} e_\lambda^x(t) cos_\lambda^{(y)}(t)
=\sum_{n=0}^\infty B_{n,\lambda}^{(c)}(x,y) \frac{t^n}{n!},
\end{split}\end{equation}
and
\begin{equation}\begin{split}\label{eq19}
\frac{t}{e_\lambda^{\frac{1}{2}}(t)-e_\lambda^{-\frac{1}{2}}(t)} e_\lambda^x(t) sin_\lambda^{(y)}(t)
=\sum_{n=0}^\infty B_{n,\lambda}^{(s)}(x,y) \frac{t^n}{n!}.
\end{split}\end{equation}
Therefore, by \eqref{eq16}, \eqref{eq17}, \eqref{eq18} and \eqref{eq19}, we obtain the following theorem.

\begin{thm}
For $n\geq 0$, we have
\begin{equation*}\begin{split}
\frac{B_{n,\lambda}(x+iy)+B_{n,\lambda}(x-iy) }{2} = B_{n,\lambda}^{(c)}(x,y),
\end{split}\end{equation*}
and
\begin{equation*}\begin{split}
\frac{B_{n,\lambda}(x+iy)-B_{n,\lambda}(x-iy) }{2i} = B_{n,\lambda}^{(s)}(x,y).
\end{split}\end{equation*}
\end{thm}

From \eqref{eq10}, \eqref{eq14} and \eqref{eq15}, we note that
\begin{equation}\begin{split}\label{eq20}
&\frac{t}{e_\lambda^{\frac{1}{2}}(t)-e_\lambda^{-\frac{1}{2}}(t)} e_\lambda^x(t) cos_\lambda^{(y)}(t)\cr
=& \sum_{l=0}^\infty B_{l,\lambda}(x) \frac{t^l}{l!} 
\sum_{m=0}^\infty \frac{(-1)^m}{(2m)!} \left( \frac{y}{\lambda}\right)^{2m} \left( \log(1+\lambda t )\right)^{2m} \cr
=&\sum_{l=0}^\infty B_{l,\lambda}(x) \frac{t^l}{l!} 
 \sum_{m=0}^\infty (-1)^m y^{2m} \lambda^{-2m}
\sum_{k=2m}^\infty S_1(k,2m) \lambda^k \frac{t^k}{k!}\cr
=& \sum_{l=0}^\infty B_{l,\lambda}(x) \frac{t^l}{l!} 
\sum_{k=0}^\infty \bigg(\sum_{m=0}^{\left[ \frac{k}{2} \right]}(-1)^m y^{2m} \lambda^{k-2m}
 S_1(k,2m)  \bigg)\frac{t^k}{k!}\cr
=&\sum_{n=0}^\infty \bigg(\sum_{k=0}^n  \sum_{m=0}^{\left[ \frac{k}{2} \right]}
\binom{n}{k} B_{n-k,\lambda}(x) (-1)^m y^{2m} \lambda^{k-2m}
 S_1(k,2m)  \bigg)\frac{t^n}{n!},
\end{split}\end{equation}
where $S_1(k,l)$ are the Stirling numbers of the first kind.
By the same method as in \eqref{eq20}, we get
\begin{equation}\begin{split}\label{eq21}
&\frac{t}{e_\lambda^{\frac{1}{2}}(t)-e_\lambda^{-\frac{1}{2}}(t)} e_\lambda^x(t) sin_\lambda^{(y)}(t)\cr
=&\sum_{l=0}^\infty B_{l,\lambda}(x) \frac{t^l}{l!}
\sum_{k=1}^\infty \bigg(\sum_{m=0}^{\left[ \frac{k-1}{2} \right]} (-1)^m y^{2m+1} \lambda^{k-2m-1}
 S_1(k,2m+1)  \bigg)\frac{t^k}{k!}\cr
=&\sum_{n=1}^\infty \bigg(\sum_{k=1}^n  \sum_{m=0}^{\left[ \frac{k-1}{2} \right]}
\binom{n}{k} B_{n-k,\lambda}(x) (-1)^m y^{2m+1} \lambda^{k-2m-1}
 S_1(k,2m+1)  \bigg)\frac{t^n}{n!}.
\end{split}\end{equation}
Therefore, by \eqref{eq18}, \eqref{eq19} , \eqref{eq20} and \eqref{eq21},
we obtain the following theorem.

\begin{thm}
For $n \in \mathbb{N} \cup \{0\}$, we have
\begin{equation*}\begin{split}
B_{n,\lambda}^{(c)}(x,y)= \sum_{k=0}^n  \sum_{m=0}^{\left[ \frac{k}{2} \right]}
\binom{n}{k} B_{n-k,\lambda}(x) (-1)^m y^{2m} \lambda^{k-2m}
 S_1(k,2m).
\end{split}\end{equation*}
In addition,
\begin{equation*}\begin{split}
 B_{0,\lambda}^{(s)}(x,y)=0,
\end{split}\end{equation*}
\begin{equation*}\begin{split}
B_{n,\lambda}^{(s)}(x,y)= \sum_{k=1}^n  \sum_{m=0}^{\left[ \frac{k-1}{2} \right]}
\binom{n}{k} B_{n-k,\lambda}(x) (-1)^m y^{2m+1} \lambda^{k-2m-1}
 S_1(k,2m+1),
\end{split}\end{equation*}
where $n$ is a positive integer.
\end{thm}

We observe that
\begin{align}\label{eq22}
\sum_{n=0}^\infty B_{n,\lambda}^{(c)}(x,0)\frac{t^n}{n!}
=& \frac{t}{e_\lambda^{\frac{1}{2}}(t)-e_\lambda^{-\frac{1}{2}}(t)} e_\lambda^x(t)\cr
=& \frac{t}{e_\lambda(t)-1} e_\lambda^{x+\frac{1}{2}}(t)\cr
=& \sum_{n=0}^\infty \beta_{n, \lambda}\left(x+\frac{1}{2}\right) \frac{t^n}{n!}.
\end{align}
Therefore, by \eqref{eq22}, we obtain the following theorem.

\begin{thm}
For $n\geq 0$, we have
\begin{equation*}\begin{split}
B_{n,\lambda}^{(c)}(x,0)= \beta_{n,\lambda}\left(x+ \frac{1}{2}\right).
\end{split}\end{equation*}
\end{thm}

From  \eqref{eq18}, we note that
\begin{equation}\begin{split}\label{eq23}
&e_\lambda^x(t) cos_\lambda^{(y)}(t)=\frac{1}{t}\left( e_\lambda^{\frac{1}{2}}(t)-e_\lambda^{-\frac{1}{2}}(t)\right)
\sum_{l=0}^\infty B_{l,\lambda}^{(c)}(x,y) \frac{t^l}{l!}\cr
=& \frac{1}{t} \sum_{n=1}^\infty \bigg( \sum_{l=0}^n \binom{n} {l}
\bigg(\left(\frac{1}{2}\right)_{n-l,\lambda}- \left(-\frac{1}{2}\right)_{n-l,\lambda} \bigg)B_{l,\lambda}(x,y)\bigg)
\frac{t^n}{n!}\cr
=&  \sum_{n=0}^\infty \left\{\frac{1}{n+1}\sum_{l=0}^{n+1} \binom{n+1} {l}
\bigg(\left(\frac{1}{2}\right)_{n+1-l,\lambda}- \left(-\frac{1}{2}\right)_{n+1-l,\lambda} \bigg)B_{l,\lambda}(x,y)
\right\}\frac{t^n}{n!}.
\end{split}\end{equation}

On the other hand,
\begin{equation}\begin{split}\label{eq24}
e_\lambda^x(t) cos_\lambda^{(y)}(t)
=&\sum_{l=0}^\infty (x)_{l,\lambda} \frac{t^l}{l!}  cos_\lambda^{(y)}(t)\cr
=&\sum_{l=0}^\infty (x)_{l,\lambda}\frac{t^l}{l!}
\sum_{m=0}^\infty \frac{(-1)^m}{(2m)!}  \left( \frac{y}{\lambda}\right)^{2m} \left( \log(1+\lambda t)\right)^{2m}\cr
=&\sum_{l=0}^\infty (x)_{l,\lambda}\frac{t^l}{l!} 
\sum_{m=0}^\infty  (-1)^m \lambda^{-2m} y^{2m} \sum_{k=2m}^\infty S_1(k,2m) \lambda^k
\frac{t^k}{k!}\cr
=&\sum_{l=0}^\infty (x)_{l,\lambda}\frac{t^l}{l!}
\sum_{k=0}^\infty \bigg(\sum_{m=0}^{\left[ \frac{k}{2} \right]}  (-1)^m \lambda^{k-2m} y^{2m} S_1(k,2m) \bigg) \frac{t^k}{k!}\cr
=& \sum_{n=0}^\infty \bigg(\sum_{k=0}^n\sum_{m=0}^{\left[ \frac{k}{2} \right]}
\binom{n}{k}(x)_{n-k,\lambda} (-1)^m \lambda^{k-2m} y^{2m} S_1(k,2m) \bigg) \frac{t^n}{n!}.
\end{split}\end{equation}
Therefore, by \eqref{eq23} and \eqref{eq24}, we obtain the following theorem.

\begin{thm}
For $n\geq 0$, we have
\begin{equation*}\begin{split}
&\frac{1}{n+1}\sum_{l=0}^{n+1} \binom{n+1} {l}
\bigg(\left( \frac{1}{2}\right)_{n+1-l,\lambda}- \left(-\frac{1}{2}\right)_{n+1-l,\lambda} \bigg)B_{l,\lambda}^{(c)}(x,y)\cr
=&\sum_{k=0}^n\sum_{m=0}^{\left[ \frac{k}{2} \right]}
\binom{n}{k}    (x)_{n-k,\lambda} (-1)^m \lambda^{k-2m} y^{2m} S_1(k,2m) .
\end{split}\end{equation*}
Furthermore, for $n \in \mathbb{N}$, we have
\begin{equation*}\begin{split}
&\frac{1}{n+1}\sum_{l=0}^{n+1} \binom{n+1} {l}
\bigg(\left( \frac{1}{2}\right)_{n+1-l,\lambda}- \left( -\frac{1}{2}\right)_{n+1-l,\lambda} \bigg)B_{l,\lambda}^{(s)}(x,y)\cr
=&\sum_{k=1}^n\sum_{m=0}^{\left[ \frac{k-1}{2} \right]}
\binom{n}{k}    (x)_{n-k,\lambda} (-1)^m \lambda^{k-2m-1} y^{2m+1} S_1(k,2m+1) .
\end{split}\end{equation*}
\end{thm}

By replacing $t$ by $\frac{1}{\lambda}\left(e^{\lambda t}-1 \right)$ in \eqref{eq18}, we get
\begin{equation}\begin{split}\label{eq25}
&\frac{1}{\lambda t}\left(e^{\lambda t}-1 \right)
\left( \frac{t}{e^{\frac{t}{2}}-e^{-\frac{t}{2}}} e^{xt} \cos yt\right)\cr
=& \sum_{k=0}^\infty B_{k,\lambda}^{(c)}(x,y)\frac{1}{k!}\left(e^{\lambda t}-1 \right)^k\lambda^{-k}\cr
=&\sum_{k=0}^\infty B_{k,\lambda}^{(c)}(x,y)\lambda^{-k}
\sum_{n=k}^\infty S_2(n,k)\lambda^n \frac{t^n}{n!}\cr
=&\sum_{n=0}^\infty\left( \sum_{k=0}^n \lambda^{n-k} B_{k,\lambda}^{(c)}(x,y)
 S_2(n,k) \right) \frac{t^n}{n!},
\end{split}\end{equation}
where $S_2(n,k)$ are the Stirling numbers of the second kind.
On the other hand,
\begin{equation}\begin{split}\label{eq26}
&\frac{1}{\lambda t}\left(e^{\lambda t}-1 \right)
\left( \frac{t}{e^{\frac{t}{2}}-e^{-\frac{t}{2}}} e^{xt} \cos yt\right)\cr
=&  \sum_{l=0}^\infty \frac{\lambda^l}{l+1}\frac{t^l}{l!} 
\sum_{m=0}^\infty \bigg(\sum_{l=0}^{\left[ \frac{m}{2} \right]}\binom{m}{2l} (-1)^l y^{2l}  B_{m-2l}(x)\bigg)\frac{t^m}{m!}\cr
=& \sum_{n=0}^\infty \bigg(\sum_{m=0}^n \frac{\lambda^{n-m}}{n-m+1}\binom{n}{m}
\sum_{l=0}^{\left[ \frac{m}{2} \right]}\binom{m}{2l}  (-1)^l y^{2l}  B_{m-2l}(x)\bigg)\frac{t^n}{n!}.
\end{split}\end{equation}
Therefore, by \eqref{eq25} and \eqref{eq26}, we obtain the following theorem.
\begin{thm}
For $n\geq 0$, we have
\begin{equation*}\begin{split}
\sum_{k=0}^n\lambda^{n-k} B_{k,\lambda}^{(c)}(x,y) S_2(n,k)
= \sum_{m=0}^n\sum_{l=0}^{\left[ \frac{m}{2} \right]} \frac{\lambda^{n-m}}{n-m+1}\binom{n}{m}
\binom{m}{2l}  (-1)^l y^{2l}  B_{m-2l}(x).
\end{split}\end{equation*}
\end{thm}
Let us replace $t$ by $\frac{1}{\lambda } \log (1+\lambda t)$ in \eqref{eq01}. Then we have
\begin{equation}\begin{split}\label{eq27}
\frac{\log(1+ \lambda t)}{\lambda t} \frac{t}{e_\lambda^{\frac{1}{2}}(t)-e_\lambda^{-\frac{1}{2}}(t)}
e^{x+iy}_\lambda (t)
=& \sum_{k=0}^\infty B_k(x+iy) \lambda^{-k} \frac{\left( \log(1+\lambda t)\right)^k}{k!}\cr
=& \sum_{k=0}^\infty B_k(x+iy) \lambda^{-k} \sum_{n=k}^\infty S_1(n,k) \lambda^n \frac{t^n}{n!}\cr
=& \sum_{n=0}^\infty \bigg( \sum_{k=0}^n \lambda^{n-k} B_k (x+iy) S_1(n,k) \bigg) \frac{t^n}{n!}.
\end{split}\end{equation}
We recall here that the Bernoulli numbers of the second are given by

\begin{equation}\label{eq27-1}
\frac{t}{\log(1+t)}=\sum_{n=0}^{\infty}b_n \frac{t^n}{n!}.
\end{equation}

Then, from \eqref{eq18}, \eqref{eq19} and \eqref{eq27}, we have
\begin{equation}\begin{split}\label{eq28}
&\sum_{n=0}^{\infty} B_{n,\lambda}^{(c)} (x,y)\frac{t^n}{n!}\cr
&= \sum_{l=0}^{\infty}b_l\lambda^l\frac{t^l}{l!}\sum_{m=0}^{\infty}\bigg(\sum_{k=0}^m \lambda^{m-k}S_1(m,k) \frac{B_k(x+iy)+B_k(x-iy)}{2}\bigg)\frac{t^m}{m!}\cr
&=\sum_{n=0}^{\infty}\bigg(\sum_{m=0}^{n}\sum_{k=0}^{m}\binom{n}{m}b_{n-m}\lambda^{n-k}
S_1(m,k)\frac{B_k(x+iy)+B_k(x-iy)}{2}\bigg)\frac{t^n}{n!},
\end{split}\end{equation}
and
\begin{equation}\begin{split}\label{eq29}
&\sum_{n=0}^{\infty} B_{n,\lambda}^{(s)} (x,y)\frac{t^n}{n!}\cr
&= \sum_{l=0}^{\infty}b_l\lambda^l\frac{t^l}{l!}\sum_{m=0}^{\infty}\sum_{k=0}^m \lambda^{m-k}S_1(m,k) \left( \frac{B_k(x+iy)-B_k(x-iy)}{2 i} \right)\frac{t^m}{m!}\cr
&=\sum_{n=0}^{\infty}\bigg(\sum_{m=0}^{n}\sum_{k=0}^{m}\binom{n}{m}b_{n-m}\lambda^{n-k}
S_1(m,k)\frac{B_k(x+iy)-B_k(x-iy)}{2 i}\bigg)\frac{t^n}{n!}.
\end{split}\end{equation}

From \eqref{eq01}, we note that
\begin{equation}\begin{split}\label{eq30}
\sum_{n=0}^\infty \left(\frac{B_n(x+iy)+B_n(x-iy)}{2}\right) \frac{t^n}{n!}
=& \frac{t}{e^{\frac{t}{2}}-e^{-\frac{t}{2}}}e^{xt} \cos yt\cr
=&\sum_{l=0}^\infty B_l(x) \frac{t^l}{l!}\sum_{m=0}^{\infty} y^{2m}(-1)^m\frac
{t^{2m}}{(2m)!}\cr
=& \sum_{n=0}^\infty \bigg(\sum_{m=0}^{\left[ \frac{n}{2}\right]}\binom{n}{2m}
B_{n-2m}(x) y^{2m} (-1)^m \bigg)\frac{t^n}{n!}.
\end{split}\end{equation}
Comparing the coefficients on both sides of \eqref{eq30}, we have
\begin{equation}\begin{split}\label{eq31}
 \frac{B_n(x+iy)+B_n(x-iy)}{2} = \sum_{m=0}^{\left[ \frac{n}{2}\right]}\binom{n}{2m}
B_{n-2m}(x) y^{2m} (-1)^m ,
\end{split}\end{equation}
where $n$ is a positive integer. By the same method as in \eqref{eq31}, we get
\begin{equation}\begin{split}\label{eq31-1}
 \frac{B_n(x+iy)-B_n(x-iy)}{2i}
 = \sum_{m=0}^{\left[ \frac{n-1}{2}\right]}\binom{n}{2m+1}
B_{n-2m-1}(x) y^{2m+1} (-1)^m ,
\end{split}\end{equation}
where $n$ is a positive integer. Therefore, by \eqref{eq28}, \eqref{eq29}, \eqref{eq31} 
\eqref{eq31-1}, we obtain the following theorem.

\begin{thm}
For $n\geq 0$, we have
\begin{equation*}\begin{split}
B_{n,\lambda}^{(c)}(x,y)
= \sum_{m=0}^{n}\sum_{k=0}^{m}\sum_{l=0}^{\left[ \frac{k}{2}\right]}\binom{n}{m}\binom{k}{2l} (-1)^l\lambda^{n-k}S_1(m,k)b_{n-m}B_{k-2l}(x) y^{2l}.
\end{split}\end{equation*}
Furthermore, for $n \in \mathbb{N}$, we have
\begin{equation*}\begin{split}
B_{n,\lambda}^{(s)}(x,y)
=\sum_{m=0}^{n}\sum_{k=0}^{m}\sum_{l=0}^{\left[ \frac{k-1}{2}\right]}\binom{n}{m}\binom{k}{2l+1}(-1)^l\lambda^{n-k}S_1(m,k)b_{n-m}B_{k-2l-1}(x) y^{2l+1}. 
\end{split}\end{equation*}
\end{thm}
For $\alpha\in \mathbb{R}$, the type $2$ degenerate Bernoulli polynomials of order $\alpha$
are defined by 
\begin{equation}\begin{split}\label{32}
\left( \frac{t}{e_\lambda^{\frac{1}{2}}(t)-e_\lambda^{-\frac{1}{2}}(t)}\right)^\alpha
e_\lambda^x (t) =\sum_{n=0}^\infty B_{n,\lambda}^{(\alpha)} (x) \frac{t^n}{n!}.
\end{split}\end{equation}
When $x=0$, $B_{n,\lambda}^{(\alpha)}=B_{n,\lambda}^{(\alpha)}(0)$ are called
the type $2$ degenerate Bernoulli numbers of order $\alpha$.
For $k \in \mathbb{N}$, let $\alpha=-k$ and $x=0$.
Then we have
\begin{equation}\begin{split}\label{eq33}
\sum_{n=0}^\infty  B_{n,\lambda}^{(-k)} \frac{t^n}{n!}
=& \frac{1}{t^k} \left( e_\lambda^{\frac{1}{2}}(t)-e_\lambda^{-\frac{1}{2}}(t)\right)^k\cr
=&  \frac{k!}{t^k}  \sum_{n=k}^\infty T_\lambda(n,k) \frac{t^n}{n!}\cr
=& \sum_{n=0}^\infty \frac{T_\lambda(n+k,k)}{\binom{n+k}{k}} \frac{t^n}{n!}.
\end{split}\end{equation}
Thus, by \eqref{eq33}, we get
\begin{equation*}\begin{split}
\binom{n+k}{k} B_{n,\lambda}^{(-k)}= T_\lambda(n+k,k),
\end{split}\end{equation*}
where $n,k$ are nonnegative integers.\\
\indent For $\alpha\in \mathbb{R}$, let us define the type $2$ degenerate cosine-Bernoulli polynomials of order $\alpha$ and the type $2$ degenerate sine-Bernoulli polynomials of order $\alpha$, repsectively by
\begin{equation}\begin{split}\label{eq36}
\bigg(\frac{t}{e_\lambda^{\frac{1}{2}}(t)-e_\lambda^{-\frac{1}{2}}(t)}\bigg)^\alpha
e_\lambda^x (t) \cos_\lambda^{(y)}(t) =\sum_{n=0}^\infty B_{n,\lambda}^{(c,\alpha)} (x,y) \frac{t^n}{n!},
\end{split}\end{equation}
and
\begin{equation}\begin{split}\label{eq37}
\bigg(\frac{t}{e_\lambda^{\frac{1}{2}}(t)-e_\lambda^{-\frac{1}{2}}(t)}\bigg)^\alpha
e_\lambda^x (t) \sin_\lambda^{(y)}(t) =\sum_{n=0}^\infty B_{n,\lambda}^{(s,\alpha)} (x,y) \frac{t^n}{n!}.
\end{split}\end{equation}
Then, we note that
\begin{equation}\begin{split}\label{eq38}
B_{n,\lambda}^{(c,\alpha)} (x,y) =\frac{B_{n,\lambda}^{(\alpha)}(x+iy)+B_{n,\lambda}^{(\alpha)}(x-iy)}{2},
\end{split}\end{equation}
where $n$ is a nonnegative integer.
\begin{equation}\begin{split}\label{eq38-1}
B_{n,\lambda}^{(s,\alpha)} (x,y) =\frac{B_{n,\lambda}^{(\alpha)}(x+iy)-B_{n,\lambda}^{(\alpha)}(x-iy)}{2i},
\end{split}\end{equation}
where $n$ is a positive integer.
Proceeding just as in \eqref{eq20} and \eqref{eq21}, we have
\begin{equation}\begin{split}\label{eq34}
&\sum_{n=0}^\infty \bigg(\frac{B_{n,\lambda}^{(\alpha)}(x+iy)+B_{n,\lambda}^{(\alpha)}(x-iy)}{2}\bigg) \frac{t^n}{n!}\cr
=& \bigg( \frac{t}{e_\lambda^{\frac{1}{2}}(t)-e_\lambda^{-\frac{1}{2}}}(t)\bigg)^\alpha
e^x_\lambda(t) \cos_\lambda^{(y)}(t)\cr
=& \sum_{n=0}^\infty \bigg(
\sum_{k=0}^n  \sum_{m=0}^{\left[ \frac{k}{2} \right]} \binom{n}{k} B_{n-k,\lambda}^{(\alpha)} (x) (-1)^m \lambda^{k-2m}y^{2m} S_1(k,2m)\bigg) \frac{ t^n}{n!},
\end{split}\end{equation}
and
\begin{equation}\begin{split}\label{eq35}
&\sum_{n=0}^\infty \bigg(\frac{B_{n,\lambda}^{(\alpha)}(x+iy)-B_{n,\lambda}^{(\alpha)}(x-iy)}{2i}\bigg) \frac{t^n}{n!}\cr
=& \bigg(\frac{t}{e_\lambda^{\frac{1}{2}}(t)-e_\lambda^{-\frac{1}{2}}(t)}\bigg)^\alpha
e^x_\lambda(t) \sin_\lambda^{(y)}(t)\cr
=& \sum_{n=1}^\infty \bigg(\sum_{k=1}^n  \sum_{m=0}^{\left[ \frac{k-1}{2} \right]} \binom{n}{k} B_{n-k,\lambda}^{(\alpha)} (x) (-1)^m \lambda^{k-2m-1}y^{2m+1} S_1(k,2m+1)\bigg) \frac{t^n}{n!}.
\end{split}\end{equation}
Therefore, by \eqref{eq38}, \eqref{eq38-1}, \eqref{eq34} and \eqref{eq35}, we obtain the following theorem.

\begin{thm}
For $n\geq 0$, we have
\begin{equation*}\begin{split}
B_{n,\lambda}^{(c,\alpha)} (x,y) 
=\sum_{k=0}^n  \sum_{m=0}^{\left[ \frac{k}{2} \right]} \binom{n}{k} B_{n-k,\lambda}^{(\alpha)} (x) (-1)^m \lambda^{k-2m}y^{2m} S_1(k,2m).
\end{split}\end{equation*}
Furthermore, for $n \in \mathbb{N}$, we have
\begin{equation*}\begin{split}
&B_{n,\lambda}^{(s,\alpha)} (x,y) \cr
=&\sum_{k=1}^n  \sum_{m=0}^{\left[ \frac{k-1}{2} \right]} \binom{n}{k} B_{n-k,\lambda}^{(\alpha)} (x) (-1)^m \lambda^{k-2m-1}y^{2m+1} S_1(k,2m+1).
\end{split}\end{equation*}
\end{thm}

For $k\in \mathbb{N}$, let $\alpha=-k$. Then, by  \eqref{eq36}, we get
\begin{equation}\begin{split}\label{eq39}
&\sum_{n=0}^\infty B_{n,\lambda}^{(c,-k)} (x,y) \frac{t^n}{n!}\cr
=& \frac{k!}{t^k}\frac{1}{k!} \left( e_\lambda^{\frac{1}{2}}(t)-e_\lambda^{-\frac{1}{2}}(t) \right)^k
e_\lambda^x(t) \cos_\lambda^{(y)}(t)\cr
=& \sum_{l=0}^\infty \frac{T_\lambda(l+k,k|x)}{\binom{l+k}{k}}  \frac{t^l}{l!}
\sum_{j=0}^\infty \bigg(\sum_{m=0}^{\left[ \frac{j}{2}\right]}(-1)^m y^{2m}\lambda^{j-2m}S_1(j,2m)\bigg) \frac{t^j}{j!}\cr
=& \sum_{n=0}^\infty \bigg( \sum_{j=0}^n  \sum_{m=0}^{\left[ \frac{j}{2}\right]}
\frac{\binom{n}{j}}{\binom{n-j+k}{k}} T_\lambda(n-j+k,k|x)(-1)^m y^{2m}\lambda^{j-2m}S_1(j,2m) \bigg) \frac{t^n}{n!} .
\end{split}\end{equation}
Therefore, by \eqref{eq39}, we obtain the following theorem.

\begin{thm}
For $k\in \mathbb{N}$ and $n \in \mathbb{N}\cup \{0\}$, we have
\begin{equation*}\begin{split}
B_{n,\lambda}^{(c,-k)}(x,y) =  \sum_{j=0}^n  \sum_{m=0}^{\left[ \frac{j}{2}\right]}
\frac{\binom{n}{j}}{\binom{n-j+k}{k}} T_\lambda(n-j+k,k|x)(-1)^m y^{2m}\lambda^{j-2m}S_1(j,2m).
\end{split}\end{equation*}
\end{thm}

From \eqref{eq11}, we define the type 2 degenerate Euler polynomials of complex variable by
\begin{equation}\begin{split}\label{eq40}
 \frac{2}{ e_\lambda^{\frac{1}{2}}(t)+e_\lambda^{-\frac{1}{2}}(t)}
e_\lambda^{x+iy}(t)=\sum_{n=0}^\infty E_{n,\lambda}(x+iy) \frac{t^n}{n!} .
\end{split}\end{equation}
From  \eqref{eq40}, we have
\begin{equation}\begin{split}\label{eq41}
\sum_{n=0}^\infty \left(\frac{E_{n,\lambda}(x+iy) + E_{n,\lambda}(x-iy)}{2}  \right) \frac{t^n}{n!}
=\frac{2e_\lambda^x(t)}{e_\lambda^{\frac{1}{2}}(t)+e_\lambda^{-\frac{1}{2}}(t) } \cos_\lambda^{(y)}(t),
\end{split}\end{equation}
and
\begin{equation}\begin{split}\label{eq42}
\sum_{n=0}^\infty \left(\frac{E_{n,\lambda}(x+iy) - E_{n,\lambda}(x-iy)}{2i}  \right) \frac{t^n}{n!}
=\frac{2e_\lambda^x(t)}{e_\lambda^{\frac{1}{2}}(t)+e_\lambda^{-\frac{1}{2}} (t)} \sin_\lambda^{(y)}(t),
\end{split}\end{equation}

Now, we define the type $2$ degenerate cosine-Euler and type $2$ degenerate sine-Euler polynomials 
as 
\begin{equation}\begin{split}\label{eq42-1}
\frac{2}{e_\lambda^{\frac{1}{2}}(t)+e_\lambda^{-\frac{1}{2}}(t)} e_\lambda^x(t) cos_\lambda^{(y)}(t)
=\sum_{n=0}^\infty E_{n,\lambda}^{(c)}(x,y) \frac{t^n}{n!},
\end{split}\end{equation}
and
\begin{equation}\begin{split}\label{eq42-2}
\frac{2}{e_\lambda^{\frac{1}{2}}(t)+e_\lambda^{-\frac{1}{2}}(t)} e_\lambda^x(t) sin_\lambda^{(y)}(t)
=\sum_{n=0}^\infty E_{n,\lambda}^{(s)}(x,y) \frac{t^n}{n!}.
\end{split}\end{equation}
By \eqref{eq11}, we see that
\begin{equation}\begin{split}\label{eq43}
&\frac{2}{e_\lambda^{\frac{1}{2}}(t)+e_\lambda^{-\frac{1}{2}}(t) }e_\lambda^x(t) \cos_\lambda^{(y)}(t)\cr
=& \sum_{l=0}^\infty E_{l,\lambda}(x) \frac{t^l}{l!}\cos_\lambda^{(y)}(t) \cr
=&\sum_{n=0}^\infty \bigg( \sum_{k=0}^n  \sum_{m=0}^{\left[ \frac{k}{2}\right]}
\binom{n}{k}E_{n-k,\lambda}(x) (-1)^m \lambda^{k-2m} y^{2m} S_1(k,2m) \bigg) \frac{t^n}{n!} ,
\end{split}\end{equation}
and
\begin{equation}\begin{split}\label{eq44}
&\frac{2}{e_\lambda^{\frac{1}{2}}(t)+e_\lambda^{-\frac{1}{2}}(t) }e_\lambda^x(t) \sin_\lambda^{(y)}(t)\cr
=&\sum_{n=0}^\infty \bigg(\sum_{k=0}^n  \sum_{m=0}^{\left[ \frac{k-1}{2}\right]}
\binom{n}{k}E_{n-k,\lambda}(x) (-1)^m \lambda^{k-2m-1} y^{2m+1} S_1(k,2m+1) \bigg) \frac{t^n}{n!}.
\end{split}\end{equation}
Therefore, by \eqref{eq42-1}, \eqref{eq42-2}, \eqref{eq43} and \eqref{eq44}, we obtain the following theorem.

\begin{thm}
For $n \in \mathbb{N}\cup \{0\}$,  we have
\begin{equation*}\begin{split}
E_{n,\lambda}^{(c)}(x,y)
=\sum_{k=0}^n  \sum_{m=0}^{\left[ \frac{k}{2}\right]}
\binom{n}{k}E_{n-k,\lambda}(x) (-1)^m \lambda^{k-2m} y^{2m} S_1(k,2m).
\end{split}\end{equation*}
Moreover, for $n \in \mathbb{N}$,
\begin{equation*}\begin{split}
E_{n,\lambda}^{(s)}(x,y)
=\sum_{k=0}^n  \sum_{m=0}^{\left[ \frac{k-1}{2}\right]}
\binom{n}{k}E_{n-k,\lambda}(x) (-1)^m \lambda^{k-2m-1} y^{2m+1} S_1(k,2m+1).
\end{split}\end{equation*}
\end{thm}

By replacing $t$ by $\frac{1}{\lambda}(e^{\lambda t} -1) $  in \eqref{eq40}, we get
\begin{equation}\begin{split}\label{eq45}
\frac{2}{e^{\frac{t}{2}}+e^{-\frac{t}{2}} }e^{(x+iy)t}
=& \sum_{k=0}^\infty E_{k,\lambda}(x+iy) \lambda^{-k}\frac{1}{k!}(e^{\lambda t} -1)^k \cr
=& \sum_{k=0}^\infty E_{k,\lambda}(x+iy) \lambda^{-k} \sum_{n=k}^\infty S_2(n,k) \lambda^n\frac{t^n}{n!}\cr
=& \sum_{n=0}^\infty \bigg( \sum_{k=0}^n E_{k,\lambda}(x+iy) S_2(n,k)
 \lambda^{n-k} \bigg) \frac{t^n}{n!}.
\end{split}\end{equation}

On the other hand,
\begin{equation}\begin{split}\label{eq46}
\frac{2}{e^{\frac{t}{2}}+e^{-\frac{t}{2}} }e^{(x+iy)t}
= \sum_{n=0}^\infty E_n (x+iy) \frac{t^n}{n!}.
\end{split}\end{equation}
Therefore, by  \eqref{eq45} and  \eqref{eq46}, we obtain the following theorem.

\begin{thm}
For $n \geq 0 $,  we have
\begin{equation*}\begin{split}
E_n(x+iy)=\sum_{k=0}^n  E_{k,\lambda}(x+iy) S_2(n,k)\lambda^{n-k}.
\end{split}\end{equation*}
\end{thm}

From  \eqref{eq46}, we can easily derive the following equation  \eqref{eq47}.
\begin{equation}\begin{split}\label{eq47}
& \sum_{n=0}^\infty \bigg(\frac{E_n(x+iy)+E_n(x-iy)}{2} \bigg) \frac{t^n}{n!} \cr
=& \frac{2}{e^{\frac{t}{2}}+ e^{-\frac{t}{2}} }e^{xt} \cos yt \cr
=& \sum_{l=0}^\infty E_l(x) \frac{t^l}{l!}
\sum_{m=0}^\infty \frac{(-1)^m y^{2m} }{(2m)!} t^{2m} \cr
=& \sum_{n=0}^\infty \bigg(\sum_{m=0}^{\left[\frac{n}{2} \right]} \binom{n}{2m}
E_{n-2m}(x) (-1)^m y^{2m} \bigg)  \frac{t^n}{n!}.
\end{split}\end{equation}
By \eqref{eq47}, we get
\begin{equation}\begin{split}\label{eq48}
\frac{E_n(x+iy)+E_n(x-iy)}{2}= \sum_{m=0}^{\left[\frac{n}{2} \right]}
\binom{n}{2m} E_{n-2m}(x) (-1)^m y^{2m},
\end{split}\end{equation}
where $n$ is a nonnegative integer.
From Theorem 2.10 and \eqref{eq48}, we have
\begin{equation}\begin{split}\label{eq49}
& \sum_{m=0}^{[\frac{n}{2}]} \binom{n}{2m} E_{n-2m}(x) (-1)^m y^{2m}\cr
=& \sum_{k=0}^n S_2(n,k) \lambda^{n-k} \left(
\frac{E_{n,\lambda}(x+iy)+E_{n,\lambda}(x-iy)}{2}\right)\cr
=& \sum_{k=0}^n S_2(n,k) \lambda^{n-k}\sum_{l=0}^k  \sum_{m=0}^{\left[\frac{l}{2} \right]}
\binom{k}{l} E_{k-l, \lambda}(x) (-1)^m \lambda^{l-2m} y^{2m}S_1(l,2m)\cr
=&  \sum_{k=0}^n \sum_{l=0}^k  \sum_{m=0}^{\left[\frac{l}{2} \right]}
S_2(n,k) \lambda^{n+l-k-2m} \binom{k}{l} E_{k-l, \lambda}(x) (-1)^m y^{2m}S_1(l,2m).
\end{split}\end{equation}
Thus, by \eqref{eq49}, we get
\begin{equation*}\begin{split}
&\sum_{m=0}^{\left[\frac{n}{2} \right]}
\binom{n}{2m} E_{n-2m}(x) (-1)^m y^{2m}\cr
=& \sum_{k=0}^n \sum_{l=0}^k  \sum_{m=0}^{\left[\frac{l}{2} \right]}
S_2(n,k) \lambda^{n+l-k-2m} \binom{k}{l} E_{k-l, \lambda}(x) (-1)^m y^{2m}S_1(l,2m).
\end{split}\end{equation*}

\medskip

\section{\bf Conclusions}
\medskip

In \cite{ref05-1}, the authors considered the degenerate Bernoulli and degenerate Euler polynomials of complex variable. By treating the real and imaginary parts separately, they were able to introduce the degenerate cosine-Bernoulli polynomials, degenerate sine-Bernoulli polynomials, degenerate cosine-Euler polynomials and  degenerate sine-Euler polynomials, and derived some interesting results for them. Actually, the degenerate Euler polynomials of complex variable are degenerate versions of the so called 'new type Euler polynomials' studied by  Masjed-Jamei, Beyki and Koepf in \cite{ref12}. Furthermore, the results in \cite{ref05-1} gave an affirmative answer to the question asked by Hac\`ene Belbachir in Mathematical Reviews (MR3808565),  "Is it possible to obtain their results by considering the classical Euler polynomials of complex variable $z$, and treating the real part and the imaginary part separately?" \\
\indent Carlitz \cite{ref01,ref02} initiated the study of degenerate versions of Bernoulli and Euler polynomials. As it turns out (see \cite{ref02-1,ref03,ref04,ref05-1,ref06,ref07,ref08,ref09} and references therein), studying degenerate versions of some special polynomials and numbers have been very fruitful and is promising. This idea of considering degenerate versions of some special polynomials is not only limited to polynomials but also can be extended to transcendental functions like gamma functions \cite{ref08}. \\
\indent In Section 2, we studied the type 2 degenerate Bernoulli and type 2 degenerate Euler polynomials of complex variable of which the latters are degenerate and type 2 versions of the aforementioned new type Euler polynomials studied in \cite{ref12}. By treating the real and imaginary parts separately, the type 2 degenerate cosine-Bernoulli and type 2 degenerate sine-Bernoulli polynomials were introduced. They were expressed in terms of the type 2 degenerate Bernoulli polynomials and Stirling numbers of the first kind. In addition, they were represented in terms of the type 2 Bernoulli polynomials and Stirling numbers of the first kind. Identities involving the type 2 degenerate cosine-polynomials (or the type 2 degenerate sine-polynomials) and Stirling numbers of the first kind were obtained. Another identity connecting  the type 2 degenerate cosine-Bernoulli polynomials, Stirling numbers of the second kind and the type 2 Bernoulli polynomials were derived. As natural extensions of the type 2 degenerate cosine-Bernoulli and type 2 degenerate sine-Bernoulli polynomials, the type 2 degenerate cosine-Bernoulli and type 2 degenerate sine-Bernoulli polynomials of order $\alpha$ were introduced. They were expressed in terms of the type 2 degenerate Bernoulli polynomials of order $\alpha$ and Stirling numbers of the second kind. In addition, the type 2 degenerate cosine-Bernoulli polynomials of negative order were represented in terms of the degenerate central factorial polynomials of the second kind and Stirling numbers of the first kind. Moreover, the type 2 degenerate cosine-Euler and type 2 degenerate sine-Euler polynomials were investigated and analogous results to  the type 2 degenerate cosine-Bernoulli and type 2 degenerate sine-Bernoulli polynomials were obtained for them.

\bigskip

{\bf Competing interests:}
The authors declare that they have no competing interests.
\bigskip

{\bf Funding:} This research received no external funding.
\bigskip

{\bf Authors' contributions:}
T.K. and D.S.K. conceived of the framework and structured the whole paper; T.K. wrote the paper; All authors read and approved the final manuscript.


\begin{thebibliography}{0}


\bibitem{ref01} L. Carlitz,
\textit{Degenerate Stirling, Bernoulli and Eulerian numbers,} Utilitas Math.  {\bf{15}} (1979), 51-88.

\bibitem{ref02} L. Carlitz,  \textit{A degenerate Staudt-Clausen theorem}, Arch. Math. (Basel) {\bf{7}} (1956), 28-33.

\bibitem{ref02-1} D. V. Dolgy, G.-W. Jang, T. Kim,  \textit{A note on degenerate central factorial polynomials of the second kind}, Adv. Stud. Contemp. Math. (Kyungshang) {\bf{29 (1)}} (2019), 7--13.


\bibitem{ref03} H. Haroon, W.A. Khan,  \textit{Degenerate Bernoulli numbers and polynomials
associated with degenerate Hermite polynomials}, Korean Math. Soc. {\bf{33 (2)}} (2018), 651-669.



\bibitem{ref04} G.-W. Jang, T. Kim,  \textit{ A note on type $2$ degenerate Euler and Bernoulli
polynomials}, Adv. Stud. Contemp. Math. (Kyungshang) {\bf{29 (1)}} (2019), 147-159.

\bibitem{ref04-1} J. Jeong, D.-J, Kang, S.-H. Rim, \textit{ Symmetry identities of Changhee polynomials of type two},
 Symmetry 2018, 10, 740. 

\bibitem{ref05} D.S. Kim, H.Y. Kim, D. Kim, T. Kim,
\textit{Identities of symmetry for type $2$ Bernoulli and Euler polynomials },
Symmetry 2019, 11, 613.

\bibitem{ref05-1} D.S. Kim, T. Kim, H. Lee,
\textit{A note on degenerate Euler and Bernoulli polynomials of complex variable },
arXiv:1908.03783 [math.NT]. https://arxiv.org/abs/1908.03783


\bibitem{ref06} T. Kim, D.S. Kim,
\textit{Degenerate central factorial numbers of the second kind},
Rev. R. Acad. Cienc. Exactas Fis. Nat. Ser. A Mat. RACSAM
(2019). https://doi.org/10.1007/s13398-019-00700-w


\bibitem{ref07} T. Kim, G.-W. Kim,
\textit{A note on degenerate gamma function and degenerate Stirling number of the second kind},
Adv. Stud. Contemp. Math. (Kyungshang)  {\bf{28 (2)}} (2018), 207-214.


\bibitem{ref08} T. Kim, D.S. Kim,
\textit{Degenerate Laplace transform and degenerate gamma function}, Russ. J. Math. Phys., {\bf{24 (2)}} (2017), 241-248.


\bibitem{ref09} T. Kim,
\textit{A note on degenerate Stirling polynomials of the second kind},
Proc. Jangjeon Math. Soc. {\bf{20 (3)}} (2017), 319-331.


\bibitem{ref10} T. Kim, D.S. Kim,
\textit{A note on type $2$ Changhee and Daehee polynomials}, Rev. R. Acad. Cienc. Exactas Fis. Nat. Ser. A Mat. RACSAM {\bf{113 (3)}} (2019), 2783-2791.



\bibitem{ref11} T. Kim, C.S. Ryoo,
\textit{Some identities for Euler and Bernoulli polynomials and their zeros}, Axioms
2018, 7, 56.



\bibitem{ref12} M. Masjed-Jamei, M.R. Beyki, W. Koepf,
\textit{A new type of Euler polynomials and numbers},
Mediterr. J. Math. {\bf{15 (3)}} (2018), Art. 138, 17pp.



\bibitem{ref13}S. Roman,
The umbral calculus. Pure and Applied Mathematics, III. Academic Press, Inc.
[Harcourt Brace Jovanovich, Publishers], New York, 1984.


\bibitem{ref14} Y. Simsek,
\textit{Identities on the Changhee numbers and Apostol-type Daehee polynomials},
Adv. Stud. Contemp. Math. (Kyungshang) {\bf{27 (2)}} (2017), 199-212.


\bibitem{ref15} W. Zhang, X. Lin,
\textit{Identities involving trigonometric functions and Bernoulli numbers},
Appl. Math. Comput. {\bf{334}} (2018), 288-294.


\end{thebibliography}
\end{document}